\documentclass{amsart}
\usepackage{amssymb,amsmath,amsthm}

\newtheorem{theorem}{Theorem}[section]
\newtheorem{lem}[theorem]{Lemma}
\newtheorem{prop}[theorem]{Proposition}

\theoremstyle{definition}
\newtheorem{definition}[theorem]{Definition}

\numberwithin{equation}{section}

\def\ord        {\mbox{\rm Ord}}
\def\tord       {\tiny{\mbox{\rm Ord}}}
\newcommand{\nn}{\nonumber}
\newcommand{\ba}{\begin{eqnarray}}
\newcommand{\ea}{\end{eqnarray}}
\newcommand{\C}{{\mathbb{C}}}
\newcommand{\R}{{\mathfrak{R}}}

\newcommand{\N}{{\mathbb{N}}}
\newcommand{\Z}{{\mathbb{Z}}}
\newcommand{\no}{\noindent}

\begin{document}
\bibliographystyle{plain}
\title[Fixed Points and Rational Functions]
{Fixed Points of Maps on the Space of Rational Functions}
\author{Edward Mosteig}
\address{Department of Mathematics, Loyola Marymount University,
Los Angeles, California 90045} \email{emosteig@lmu.edu}

\subjclass{Primary 33}

\date{\today}


\begin{abstract}
Given integers $s,t$, define a function $\phi_{s,t}$ on the space of
all formal series expansions by $\phi_{s,t} (\sum a_n x^n) = \sum
a_{sn+t} x^n$. For each function $\phi_{s,t}$, we determine the
collection of all rational functions whose Taylor expansions at zero
are fixed by $\phi_{s,t}$.  This collection can be described as a
 subspace of rational functions whose basis elements
correspond to certain $s$-cyclotomic cosets associated with the pair
$(s,t)$.
\end{abstract}

\maketitle


\section{Introduction} \label{intro}
\setcounter{equation}{0}

Let $\R$ denote the space of rational functions with complex
coefficients. The Taylor expansion at $x =0$ of $R \in \R$ can be
written as a Laurent series, i.e., \ba\label{laurent} R(x) & = &
\sum_{n \gg -\infty} a_{n}x^{n} \ea \no where $n \gg -\infty$
denotes the fact that the coefficients vanish for large negative
$n$.
 For $s, t \in \Z$, define the map $\phi_{s,t}: \R
\to \R$ by \ba \phi_{s,t}(\sum a_nx^n) = \sum a_{sn+t}x^n. \ea
Denote the standard $s$-th root of unity throughout this paper by
$\omega_s = e^{2\pi i/s}.$
 When $s$ is positive, consider the restriction  $\phi_{s,t}: \R \to \R$.
One can rewrite this map explicitly without the use of series
expansions: \ba \label{phist}\phi_{s,t}(R(x))=
\left(\frac{1}{s}\right) x^{-t/s} \sum_{j=0}^{s-1} \omega_s^{-jt}
R(\omega_s^j x^{1/s}). \ea Indeed, if $R(x) = \sum a_n x^n$, then
$R(\omega_s^j x^{1/s}) = \sum a_n \omega_s^{jn} x^{n/s}$, and so the
coefficient of $x^{(sn+t)/s}$ in the summation $\sum a_n
\omega_s^{jn} x^{n/s}$ is $a_{sn+t}\omega_s^{j(sn+t)}$. Therefore,
the coefficient of $x^n$ in $\left(\frac{1}{s}\right) x^{-t/s}
\sum_{j=0}^{s-1} \omega_s^{-jt} R(\omega_s^j x^{1/s})$ is \break
$\left(\frac{1}{s}\right) \sum_{j=0}^{s-1} \omega_s^{-jt}
a_{sn+t}\omega_s^{j(sn+t)} = \left(\frac{1}{s}\right)
\sum_{j=0}^{s-1} a_{sn+t} = a_{sn+t}$.

The map $\phi_{2,1}$ can be used in a general procedure for the
exact integration of rational functions, as described in
\cite{bomolan2}. Dynamical properties of $\phi_{2,1}$, including
kernels of the iterates, dynamics of subclasses of rational
functions, and
 fixed points are discussed in \cite{bolimomost}.
The purpose of this paper is to generalize one of the results in
\cite{bolimomost} by classifying, for each pair of integers $s,t$,
the collection of all rational functions that are fixed by
$\phi_{s,t}$.  If $s$ is an integer such that $s \le 1$,  then 0
is the only rational function fixed by $\phi_{s,t}$, unless, of
course, $(s,t)=(1,0)$, in which case $\phi_{s,t}$ is the identity.
When $s \ge 2$, however, the story is much more interesting.


\section{Cyclotomic Cosets} \label{genfun}
\setcounter{equation}{0}

In this section, we assume throughout that $s \ge 2$, $0 \le t \le
s-2$, and $R \in \R$ such that \ba \phi_{s,t}(R(x)) = R(x). \ea
Given these restrictions on $s$ and $t$, it follows that $|t/(s-1)|
< 1$.  Thus, if $n\le -1$, then $n < -t/(s-1)$, and so $sn+t< n$.
Assuming that $R(x)$ is fixed by $\phi_{s,t}$, we have that
$a_{sn+t} = a_n$ for all $n$.  Thus, if $a_n$ is nonzero for any
negative value of $n$, then there are infinitely many nonzero
coefficients of negative powers of $x$, contradicting the assumption
that $R(x)$ is of the form given in equation (\ref{laurent}).

We write $R$ in the form \ba R(x) = \sum_{n=0}^\infty f(n) x^n \ea
to emphasize the fact that the coefficients can be interpreted as
the images of a generating function $f: \N \to \C$. Since $R(x)$
is fixed by $\phi_{s,t}$, it follows that \ba f(n) = f(sn+t) \ea
for all integers $n$. The following result, which was proven on
page 202 of \cite{stan1}, elucidates the relationship between the
generating function $f$ of the coefficients of the Taylor
expansion of $R(x)$ and the representation of $R(x)$ as a quotient
of polynomials.

\begin{lem}\label{stanley}
Let $q_1, q_2, \dots, q_d$ be a fixed sequence of complex numbers,
$d \ge 1$, and $q_d \not= 0$. The following conditions on a
function $f: \N \to \C$ are equivalent:
\begin{enumerate}
\item $\sum_{n \ge 0} f(n)x^n = \frac{P(x)}{Q(x)}$ where, $Q(x) =
1 + q_1 x + q_2 x^2 + q_3 x^3 + \cdots + q_d x^d$. \item For  $n
\gg 0$,
$$f(n) = \sum_{i=1}^J P_i(n) \lambda_i^n,$$
where $1+ q_1 x + q_2 x^2 + q_3 x^3 + \cdots + q_d x^d=
\prod_{i=1}^J (1 - \lambda_i x)^{d_i}$, the $\lambda_i$'s are
distinct, and $P_i(n)$ is a polynomial in $n$ of degree less than
$d_i$.
\end{enumerate}
\end{lem}

In this section, we construct  a collection of rational functions
that are fixed by $\phi_{s,t}$,  and in the next section we use
the above lemma to justify that this collection spans the subspace
of $\R$ consisting of all rational functions that are fixed by
$\phi_{s,t}$.

The description of all the fixed points of $\phi_{s,t}$ requires
the notion of {\em cyclotomic cosets}:  given $n, \, r \in
\mathbb{N}$ with  $r \ge 1$ such that $r$ and $s$ are relatively
prime,
\ba
C_{s,r,n} & = & \{ s^{i}n \, \text{ mod } \, r:
\; i \in \mathbb{Z} \}
\ea
\no
is a finite set called the $s$-cyclotomic coset of $n$
mod $r$.  We will
characterize the fixed points $\phi_{s,t}$ using cyclotomic cosets
with a special property.  To describe this property, first define
\ba \beta_{s,t}(k) = t \left( \frac{s^k-1}{s-1} \right) \ea for
which we have the following recursive formula: \ba
\label{recursivebeta} \beta_{s,t}(j+1) = s\beta_{s,t}(j) + t. \ea

\begin{definition}
A positive integer $r$ is called {\em distinguished} with respect
to the pair $(s,t)$ if $r$ and $s$ are relatively prime and \ba r
\mid \beta_{s,t}(\mbox{Ord}(s;r)), \ea where Ord$(s;r)$ represents
the smallest positive integer $i$ such that $s^i \equiv 1 \mod r$.
We say $r=0$ is distinguished with respect to $(s,t)$ if and only
if $t=0$.  We denote the set of integers distinguished with
respect to $(s,t)$ by $\Omega(s,t)$.
\end{definition}

\begin{prop}\label{prop:omegainfinite}
For each pair $(s,t)$, the set $\Omega(s,t)$ is infinite.
\end{prop}

\begin{proof}
Since $\Omega(s,t) \subset \Omega(s,1)$, we need only show that
$\Omega(s,1)$ is infinite.  Let $r$ be a positive integer such that
$\gcd(r,s(s-1)) = 1$.  If $\alpha = \ord(s;r)$, then $s^ \alpha
\equiv 1 \mod r$; that is, $r \mid s ^ {\tord(s;r)}-1$. Since $s ^
{\tord(s;r)}-1$ is a multiple of $(s-1)$, and $r$ is relatively
prime to $(s-1)$, it follows that $r(s-1)  \mid s ^ {\tord(s;r)}-1$.
Thus, $r \mid \frac{s ^ {\tord(s;r)}-1}{s-1} =
\beta_{(s,1)}(\ord(s;r))$, and so $r$ is distinguished with respect
to $(s,1)$.
\end{proof}

For example,  consider
\[ \Omega(3,1) = \{1, 4, 5, 7, 10, 11, 13, 14,
 17, 19, 20, 23, 25, 28,
29, 31, 34, 35, 37, 38,  \dots \}.
\]
From Proposition \ref{prop:omegainfinite}, we see that $\Omega(3,1)$
contains the arithmetic sequences $\{6n+1\}$ and $\{6n+5\}$.  With a
little more effort, one can show that $\Omega(3,1)$ also contains
the arithmetic sequences $\{24n+4\}$, $\{24n+10\}$, $\{24n+14\}$,
and $\{24n+20\}$.  The smallest integer in $\Omega(3,1)$ not
contained in any of these sequences is 40. Moreover, a calculation
shows that $96n+40$, for $0 \le n \le 5$ is in $\Omega(3,1)$, but
$616= 96 \cdot 6 + 40$ is not in $\Omega(3,1)$. An interesting
question of further study is whether the sets $\Omega(s,t)$ have a
nice characterization.  For example, we might ask whether they can
be written as a (possibly infinite) union of arithmetic sequences,
as is the case for $\Omega(2,1)$, which consists precisely of all
odd natural numbers.   However,  the example $\Omega(3,1)$ suggests
that this may not be the case in general.

A generating set for the collection of fixed points of $\phi_{s,t}$
will be indexed by $s$-cyclotomic cosets $C_{s,r,n}$ where $r$ is
distinguished with respect to $(s,t)$. Note that by computing
\ba\label{eq:psiimage} \phi_{s,t} \left( \frac{1}{1 - \lambda x}
\right) & = & \frac{\lambda^t}{1 - \lambda^{{s}}x} \ea we acquire
the following formula for the iterates of $\phi_{s,t}$: \ba
\label{eq:psikimage} \phi_{s,t}^{(k)} \left( \frac{1}{1 - \lambda x}
\right) & = & \frac{\lambda^{ \beta_{s,t}(k)  }}{1 -
\lambda^{{s^k}}x}  . \ea
 For $r \ge 1$ and
$n\in\N$, define \ba \label{eq:psi} \psi_{s,t,r,n}(x)  =
 \sum_{j=1}^{\mbox{\tiny{Ord}}(s;r)}
\frac{\omega_r^{n{\beta_{s,t}(j)}}}{1 - \omega_r^{ns^j} x}
 =
 \sum_{j=1}^{\mbox{\tiny{Ord}}(s;r)}
\phi_{s,t}^{(j)} \left(\frac{1}{1-\omega_r^nx} \right) . \ea Note
that if $n=0$, then $\psi_{s,t,r,0} = 1/(1-x)$.
 If
$t=0$, then $r=0$ is distinguished with respect to $(s,t)$, and we
define \ba \psi_{s,0,0,n}(x) = 1. \ea

\begin{prop} \label{basis}
If $r$ is distinguished with respect to $(s,t)$, then
$\psi_{s,t,r,n}(x)$ is fixed by $\phi_{s,t}$.
\end{prop}

\begin{proof}
If $r > 1$ is distinguished with respect to $(s,t)$, then \ba
\phi_{s,t}^{(\mbox{\tiny Ord}(s;r)+1)} \left(
\frac{1}{1-\omega_r^nx}\right) = \phi_{s,t}
\left(\phi_{s,t}^{(\mbox{\tiny Ord}(s;r))} \left(
\frac{1}{1-\omega_r^nx}\right)\right) = \phi_{s,t}^{} \left(
\frac{1}{1-\omega_r^nx}\right) ,\nn \ea and so \ba
\phi_{s,t}(\psi_{s,t,r,n}(x)) = \phi_{s,t} \left(
\sum_{j=1}^{\mbox{\tiny{Ord}}(s;r)} \phi_{s,t}^{(j)}
\left(\frac{1}{1-\omega_r^nx} \right)  \right) = \left(
\sum_{j=1}^{\mbox{\tiny{Ord}}(s;r)} \phi_{s,t}^{(j+1)}
\left(\frac{1}{1-\omega_r^nx} \right)  \right) = \psi_{s,t,r,n}(x).
\nn \ea \no Thus $\psi_{s,t,r,n}(x)$ is fixed by $\phi_{s,t}$. Since
constants are fixed by $\phi_{s,0}$, it follows that
$\psi_{s,0,0,n}$ is fixed by $\phi_{s,0}$. Since $r=0$ is
distinguished only with respect to $t=0$, we have shown the result
holds in all possible cases.
\end{proof}


\section{The Space of Fixed Points of $\phi_{s,t}$}\label{fixedpoints}

We now classify all the fixed points of $\phi_{s,t}$ for all
integers $s,t$.   To do so, we first demonstrate a bijective
correspondence between fixed points of $\phi_{s,t}$ and
$\phi_{s,t+u(s-1)}$ where $u$ is an arbitrary integer.

\begin{lem} \label{correspondence}
For all integers $s,t,u$, the rational function $R(x)$ is a fixed
point of $\phi_{s,t}$ iff $x^{-u}R(x)$ is a fixed point of
$\phi_{s,t+u(s-1)}$.
\end{lem}

\begin{proof}
Using equation (\ref{phist}), one can show directly that for any
integers $s,t,u$,
 \ba
\phi_{s,t}(R(x)) = x^u \phi_{s,t+(s-1)u}(x^{-u}R(x)),\nn \ea and so
\begin{eqnarray*}
\phi_{s,t}(R(x)) = R(x) &  \Leftrightarrow &
x^{u} \phi_{s,t+(s-1)u} (x^{-u}R(x)) = R(x) \\
& \Leftrightarrow & \phi_{s,t+(s-1)u} (x^{-u}R(x)) = x^{-u}R(x) .
\nn
\end{eqnarray*}\no
\end{proof}

Given this correspondence, we only have to compute the fixed points
of $\phi_{s,t}$ in case $0 \le t \le s-2$. Once this is
accomplished, to compute the fixed points of $\phi_{s,t}$ for
arbitrary $t$, we only need to find $t',u$ such that $0 \le t' \le
s-2$ and $t = t' + u(s-1)$, and then use the correspondence. The
following result provides the missing component of this scheme, thus
allowing us to compute the fixed points $\phi_{s,t}$ for any
integers $s$ and $t$.

\begin{prop}
\label{fixedpt} Suppose $s \ge 2$ and $0 \le t \le s-2$. A rational
function is fixed by $\phi_{s,t}$ if and only if it is a linear
combination of the functions $\psi_{s,t,r,n}(x)$ where $r$ is
distinguished with respect to $(s,t)$ and $n$ is relatively prime to
$r$.
\end{prop}

\begin{proof}
We showed in Proposition \ref{basis}  that if $r$ is distinguished
with respect to $(s,t)$, then $\psi_{s,t,r,n}(x)$ is fixed by
$\phi_{s,t}$, and so every linear combination of such functions
must be fixed by $\phi_{s,t}$.

To prove the converse, we consider a rational function $R(x)$ fixed
by $\phi_{s,t}$, and express it as \ba R(x) = C(x) +
\frac{P(x)}{Q(x)}\ea where $C(x),P(x),Q(x)$ are polynomials such
that $P(x)$ and $Q(x)$ are relatively prime with $\deg P(x) < \deg
Q(x)$.   Our first goal is to show that the poles of $R(x)$ must be
simple. We write
 \ba \frac{P(x)}{Q(x)}
 = \sum_{n=0}^\infty f(n) x^n \ea where  $f(n)$ is the generating function for
$P(x)/Q(x)$.   Since $f(n) = f(sn+t)$, we have by Lemma
\ref{stanley}, $f(sn+t) = \sum
P_i(sn+t)\lambda_i^t(\lambda_i^s)^n$ and\ba Q(x) = \prod_{i=1}^J
(1 - \lambda_i x)^{d_i} = \prod_{i=1}^J (1- \lambda_i^sx)^{e_i},
\ea \no and so \ba \label{lambdaispower} \{ \lambda_1, \dots,
\lambda_J \} = \{ \lambda_1^s, \dots, \lambda_J^s \}. \ea \no Thus
the set $\{ \lambda_{1}, \cdots, \lambda_{J} \}$ is permuted by
the
 map $z\mapsto z^s$, and so each $\lambda_j$ is a primitive $r_j$-th root
 of unity where $r_j$ is a positive integer.
Moreover, since $\{ \lambda_1, \dots, \lambda_J\}$ is permuted by
the map $z \mapsto z^s$, it follows that for each $1 \le j \le J$,
there exists a positive integer $\ell$ such that $\lambda_j^{s\ell}
= \lambda_j$ (after applying the map $z \mapsto z^s$ multiple
times). Therefore, $\lambda_j^{s\ell-1} = 1$, and so $r_j \mid
s\ell-1$. Thus $r_j$ and $s$ are relatively prime.

Let $M = $ lcm$(r_1, \dots, r_j)$ and for $a\in \N$, define
$$R_a = \{ m \in \N : m
\equiv a \; \text{ mod } M \}.$$ Let $f_a = f \Big{|}_{R_{a}}$ be
the restriction of the function $f : \N \to \C$ to the set $R_a$.
Then
$$f_a(a+jM) = \sum_{i=1}^J P_i(a+jM) \lambda_i^{a+jM} =
\sum_{i=1}^J P_i(a+jM) \lambda_i^{a},$$ and so each $f_a$ has a
representation as a polynomial in the variable $j$ since
$\lambda_i^a$ is constant on the set $R_a$. We denote the natural
extension of this map to an element of the polynomial ring $\C[j]$
by $F_a$. Note that the restriction of $F_a$ to $\N$ need not be
$f$ in general. Our goal is to prove that each $F_a$ is a constant
function, with corresponding constant denoted by $c_a$.
 Once this is shown, we have
\ba \frac{P(x)}{Q(x)} = \sum_{n=0}^\infty f(n) x^n =
 \sum_{a=0}^{M-1} c_a \sum_{j=0}^{\infty}
 x^{a+jM} = \sum_{a=0}^{M-1} \frac{c_ax^a}{1-x^M},
\ea and so $P(x)/Q(x)$ is a rational function with only simple
poles, as desired.

It remains to show  that each polynomial map $F_a: \C \to \C$ is a
constant function.   For each positive integer $n$, define \ba S_n =
\{ \beta_{s,t}^{(j)}(n) : j \in \N  \}. \ea \no We say that $a$ has
an {\em infinite cross-section} if $R_a \cap S_n$ is an infinite set
for some $n \in \N$. We proceed by considering two cases, depending
on whether $a$ has an infinite cross-section or not \\

\no {\bf Case 1}:  Suppose $a$ has an infinite cross-section,
i.e., $R_a \cap S_n$ is an infinite set.  Since $f(j) = f(sj+t)$
for all $j\in \N$ , $F_a$ is constant on $R_a \cap S_n$.  Since
$R_a \cap S_n$
is an infinite set, $F_a$ is a constant polynomial. \\

\no {\bf Case 2}:  Suppose $a$ does not have an infinite
cross-section, i.e.,  $R_a \cap S_n$ is finite for all positive
integers $n$. Then $R_a \cap S_n$ must be nonempty for infinitely
many values of $n$.  Since there are only finitely many distinct
sets of the form $R_b$, it follows that for  each $S_n$, there
exists $b \in \N$ such that $R_b \cap S_n$ is infinite.  Moreover,
since there are only finitely many choices for $R_b$, there is at
least one $b \in \N$ such  that there exist infinitely many values
of $n$
 where $R_a \cap S_n$ is nonempty and $R_b \cap S_n$ is
infinite.  Since $b$ has an infinite cross-section, an application
of Case 1 demonstrates that the restriction of $f$ to $R_b$ is the
constant function $c_b$. Since $f$ is constant on each $S_n$, the
restriction of $f$ to $S_n$ is the constant $c_b$. Thus $F_a$
achieves the value $c_b$ infinitely many times, and so $F_a$ must
be a constant polynomial.

Thus in either case, we have that $F_a$ is a constant polynomial,
and so the poles of $R$ must be simple.
 Using this fact, we can decompose $R(x)$ using partial fractions: \ba
\label{partialfrac} R(x) = C(x) + \sum_{j=1}^J \frac{\alpha_j}{1-
\lambda_jx}. \ea Via (\ref{eq:psiimage}), an application of
$\phi_{s,t}$ yields \ba R(x) = \phi_{s,t}(R(x)) = \phi_{s,t}(C(x)) +
\sum_{j=1}^J \frac{\alpha_j \lambda_j^t}{1- \lambda_j^s x}. \ea Each
rational function has a unique decomposition, and since $\phi_{s,t}$
maps polynomials to polynomials, \ba C(x) = \phi_{s,t}(C(x)) \ea and
\ba \label{Rreduced} \sum_{j=1}^J \frac{\alpha_j}{1- \lambda_jx} =
\sum_{j=1}^J \frac{\alpha_j \lambda_j^t}{1- \lambda_j^s x} =
\phi_{s,t}\left( \sum_{j=1}^J \frac{\alpha_j}{1- \lambda_jx}
\right). \ea If $t>0$, it is easy to see that no nonzero polynomial
is fixed by $\phi_{s,t}$, in which case $C(x) =0$.  If $t=0$, then
the only polynomials fixed by $\phi_{s,t}$ are constant, and so
$C(x)$ is a constant multiple of $\psi_{s,0,0,n} = 1$.

Now we only have left to show that the second summand in
(\ref{partialfrac}) is a linear combination of functions of the form
$\psi_{s,t,r,n}$. To do this, we begin by showing that each $r_k$ is
distinguished with respect to $(s,t)$.  We have already shown that
$r_k$ and $s$ are relatively prime for each $k$.   Using
(\ref{eq:psikimage}), multiple iterations of $\phi_{s,t}$ to
(\ref{Rreduced}) yield \ba \sum_{j=1}^J \frac{\alpha_j}{1-
\lambda_jx}  =
\phi_{s,t}^{(\mbox{\tiny{Ord}}(s;r_k))}\left(\sum_{j=1}^J
\frac{\alpha_j}{1- \lambda_jx}  \right)
 = \sum_{j=1}^J
\frac{\alpha_j\lambda_j^{\beta_{s,t}(\mbox{\tiny{Ord}}(s;r_k))}}{1-
\lambda_j^{s^{\mbox{\tiny{Ord}}(s;r_k)}}x}.
\ea

The term corresponding to $j=k$ in the first of these three
expressions is \ba \frac{\alpha_k}{1-\lambda_kx}, \ea and the
corresponding term in the last of these three expressions is \ba
\frac{\alpha_k\lambda_k^{\beta_{s,t}(\mbox{\tiny{Ord}}(s;r_k))}}{1-
\lambda_k^{s^{\mbox{\tiny{Ord}}(s;r_k)}}x} =
\frac{\alpha_k\lambda_k^{\beta_{s,t}(\mbox{\tiny{Ord}}(s;r_k))}}{1-
\lambda_kx}. \ea Thus \ba
\lambda_k^{\beta_{s,t}(\mbox{\tiny{Ord}}(s;r_k))} =1. \ea
Therefore, $r_k \mid \beta_{s,t}(\mbox{{Ord}}(s;r_k))$, and so
$r_k$ is distinguished with respect to $(s,t)$.

Now that we've shown that each $r_k$ is distinguished with respect
to $(s,t)$, group terms in the sum \ba \label{unique} \sum_{j=1}^J
\frac{\alpha_j}{1- \lambda_jx} \ea according to the orbits of the
map $z \mapsto z^s$ on the set $\{\lambda_1, \dots, \lambda_J \}$.
Since $r_k \mid \beta_{s,t}(\mbox{{Ord}}(s;r_k))$ for each $k$, we
know that the sum of terms in (\ref{unique}) corresponding to a
single orbit must be of the form \ba {\mathcal O}(k) = \sum_{i=1}^m
\phi_{s,t}^{(i)} \left(\frac{\alpha_k}{1 - \lambda_k x}   \right),
\ea where $m$ is the length of the orbit of $\lambda_k$ under the
map $z \mapsto z^s$. That is, $m$ is the smallest positive integer
such that $\lambda_k^{s^m} = 1$, and so $m = \mbox{{Ord}}(s;r_k)$.
Moreover, $\lambda_k$ is a primitive $r_k$-th root of unity, and so
it must be of the form $\lambda_k = (\omega_{r_k})^n$ for some $n\in
\N$ such that $r_k$ and $n$ are relatively prime. Thus \ba {\mathcal
O}(k) = \alpha_k \left(\sum_{i=1}^{\mbox{\tiny Ord}(s;r_k)}
\phi_{s,t}^{(i)} \left(\frac{1}{1 - \omega_{r_k}^n x} \right)
  \right) = \alpha_k \psi_{s,t,r_k,n}(x),
\ea and so (\ref{unique}), and hence (\ref{partialfrac}), is a
linear combination of rational functions of the form
$\psi_{s,t,r,n}$.
\end{proof}

It turns out that the collection of rational functions of the form
$\psi_{s,t,r,n}$ does not form a basis of fixed points.  The lemma
below shows that there is redundancy in the collection. Since
cyclotomic cosets have many different representations, we must
compare the ways in which $\psi_{s,t,r,n}$ and $\psi_{s,t,r,n'}$ are
defined for two distinct representations $C_{s,r,n}$ and
$C_{s,r,n'}$ of the same coset. Although we have not defined
$\psi_{s,t,r,n}$ to be invariant with respect to different
representations, they will be the same up a constant multiple.

\begin{lem}\label{scalarcyclotomic}
If $C_{s,r,n} = C_{s,r,n'}$, then $\psi_{s,t,r,n}$ and
$\psi_{s,t,r,n'}$ are scalar multiples of one another.
\end{lem}

\begin{proof}
With the aid of (\ref{recursivebeta}), we compute
\begin{eqnarray*}
\psi_{s,t,r,ns}(x)
& = &  \sum_{j=1}^{\mbox{\tiny Ord}(s;r)} \frac{\omega_r^{ns\beta_{s,t}(j)}}{1-\omega_r^{sns^j}x} \\
& = &  \sum_{j=1}^{\mbox{\tiny Ord}(s;r)} \frac{\omega_r^{n(\beta_{s,t}(j+1)-t)}}{1-\omega_r^{ns^{j+1}}x} \\
& = &  \omega_r^{-nt} \sum_{j=1}^{\mbox{\tiny Ord}(s;r)} \frac{\omega_r^{n\beta_{s,t}(j+1)}}{1-\omega_r^{ns^{j+1}}x} \\
& = &  \omega_r^{-nt} \sum_{j=1}^{\mbox{\tiny Ord}(s;r)} \phi_{s,t}^{(j+1)} \left( \frac{1}{1-\omega_r^{n}x} \right) \\
& = &  \omega_r^{-nt} \phi_{s,t}(\psi_{s,t,r,n}(x)) \\
& = & \omega_r^{-nt} \psi_{s,t,r,n}(x).
\end{eqnarray*}
\no Thus $\psi_{s,t,r,s^in}$ and $\psi_{s,t,r,n}$ are scalar
multiples of one another for all $i \in \N$.  If $C_{s,r,n} =
C_{s,r,n'}$, then for some $i \in \N$, we have $n' \equiv s^in \mod
r$. By this equivalence, $\psi_{s,t,r,n'} = \psi_{s,t,r,s^in}$, and
so the result follows.
\end{proof}

Using  Lemma \ref{scalarcyclotomic}, we can show that if two of
functions of the form $\psi_{s,t,r,n}$ have a pole in common, then
they are actually the same up to a scalar multiple.  The following
lemma leads us this result.

\begin{lem} \label{commonpoles}
Suppose $r_i$ is a positive integer that is distinguished with
respect to $(s,t)$, and $n_i$ is a positive integer relatively prime
to $r_i$ for $i=1,2$. If  $\psi_{s,t,r_1,n_1}$ and
$\psi_{s,t,r_2,n_2}$ have a pole in common, then $r_1=r_2$ and
$C_{s,r_1,n_1} = C_{s,r_2,n_2}$.
\end{lem}

\begin{proof}
Note that $\psi_{s,t,r,n}$ has poles at $\omega_r^{-ns^j}$ for $0
\le j \le \ord(s;r)$; that is, $\psi_{s,t,r,n}$ has poles at
$\omega_r^{-c}$ where $c \in C_{s,r,n}$. Suppose
$\psi_{s,t,r_1,n_1}$ and $\psi_{s,t,r_2,n_2}$ have a pole in common;
that is, $e^{-2\pi i c_1/r_1} = e^{-2\pi i c_2/r_2}$, where $c_i \in
C_{s,r_i,n_i}$. Thus, $c_1/r_1  - c_2/r_2 \in \Z$. Without loss of
generality, we can choose $1 \le c_i < r_i$, in which case $ 0 <
c_1/r_1 < 1$, and so $c_1/r_1 = c_2/r_2$.
 Since
$\gcd(r_i,n_i) = 1$ and $c_i = s^{j_i} n_i \mod r_i$ for some $j_i
\in \N$, it follows that $c_i$ and $r_i$ are relatively prime, and
so $c_1 = c_2$ and $r_1=r_2$. Therefore, $ s^{j_1}n_1 = s^{j_2}n_2
\mod r$ (where $r=r_1=r_2$), and so $C_{s,r,n_1} = C_{s,r,n_2}$.
\end{proof}

We now precisely describe the redundancy in the collection $\{
\psi_{s,t,r,n}\}$ for fixed $s$ and $t$.  We begin by defining an
equivalence relation $\sim_{s,r}$ on $(C_{s,r,n}-\{0\})$ by $n_1
\sim_{s,r} n_2$ if $C_{s,r,n_1} = C_{s,r,n_2}$. Let $\Lambda_{s,r}$
be a collection of coset representatives (all chosen to be less than
$r$) of $(C_{s,r,n}-\{0\})/\sim_{s,r}$.  That is, $\Lambda_{s,r}$ is
maximal set consisting of positive integers such that no two are in
the same cyclotomic coset.

\begin{theorem}
 Suppose $s \ge 2$ and $0 \le t \le s-2$.
The function $1/(1-x)$ together with the collection of all
$\psi_{s,t,r,n}$ where $r$ is distinguished with respect to $(s,t)$
and $n \in \Lambda_{s,r}$ form a basis for the set of all rational
functions that are fixed points of $\phi_{s,t}$.
\end{theorem}

\begin{proof}
The case  $n=0$ corresponds to the function $1/(1-x)$.  We now
consider the case $n>0$.  Given an integer $r$ that is distinguished
with respect to $(s,t)$, and an integer $n$ that is relatively prime
to $r$, there exists $n' \in \Lambda_{s,r}$ such that $C_{s,r,n} =
C_{s,r,n'}$, in which case by Lemma \ref{scalarcyclotomic},
$\psi_{s,t,r,n}$ and $\psi_{s,t,r,n'}$ are scalar multiples of one
another.  Thus, by Proposition \ref{fixedpt}, this collection spans
the space of rational functions fixed by $\phi_{s,t}$.

Suppose $\psi_{s,t,r_1,n_1}$ and $\psi_{s,t,r_2,n_2}$ have a pole in
common where $n_i \in \Lambda_{s,r_i}$.  Then by Lemma
\ref{commonpoles}, $r_1=r_2$ and $C_{s,r_1,n_1} = C_{s,r_2,n_2}$.
Thus by the definition of $\Lambda_{s,r_1} = \Lambda_{s,r_2}$,
$n_1=n_2$.  Therefore, none of the elements of the collection have a
pole in common, and so no nontrivial linear combination of elements
of this collection can be zero.
\end{proof}


\begin{thebibliography}{99}

\bibitem{bolimomost}
G. Boros, J. Little, V. Moll, E. Mosteig, R. Stanley, A map on the
space of rational functions, {\em Rocky Mountain Journal of
Mathematics}, to appear.


\bibitem{bomolan2}
G. Boros, V. Moll, Landen transformations and the integration of
rational functions, {\em Math. Comp.} {\bf 71} (2002), 649-668.

\bibitem{mollnot}
V. Moll, The evaluation of integrals: a personal story, {\em
Notices AMS}  {\bf 49} (2002), 311-317.

\bibitem{stan1}
R. Stanley, {\em Enumerative Combinatorics, Volume} 1,  Cambridge
Studies in Advanced Mathematics, Cambridge University Press, 1997.




\end{thebibliography}
\end{document}